\documentstyle{amsart}
\def\dj{d\kern-.30em\raise1.25ex\vbox{\hrule width .3em height .03em}}
\def\Dj{D\kern-.70em\raise0.75ex\vbox{\hrule width .3em height .03em}
\kern.03em}
\newcommand{\e}{\epsilon}
\newcommand{\k}{\kappa}
\newcommand{\id}{\mbox{\shape{n}\selectfont id}}
\newtheorem{lem}{Lemma}
\newtheorem{pro}[lem]{Proposition}
\theoremstyle{definition}
\newtheorem{defn}{Definition}
\begin{document}
\title{On braided quantum groups}
\author{Mi\'co \Dj ur\Dj evi\'c}
\address{Instituto de Matematicas, UNAM, Area de la Investigacion
Cientifica, Circuito Exterior, Ciudad Universitaria, M\'exico DF,
CP 04510, M\'EXICO\newline
\indent {\it Written In}\newline
\indent Faculty of Physics, University of Belgrade,
Belgrade, SERBIA\newline
\indent {\it Final Version}\newline
\indent Centro de Investigaciones Teoricas, Facultad de Estudios Superiores
Cuautitlan, Cuautitlan Izcalli, M\'EXICO}
\email{miko@@unamvm1.dgsca.unam.mx}
\maketitle
\begin{abstract}
A braided generalization of the concept of Hopf algebra (quantum group)
is presented. The generalization overcomes an inherent geometrical
inhomogeneity of quantum groups, in the sense of allowing completely
pointless objects. All braid-type equations appear as
a consequence of initial axioms.
Braided counterparts of basic algebraic relations
between fundamental
entities of the standard theory are found.
\end{abstract}
%\tableofcontents
\section{Introduction}
{\renewcommand{\thepage}{}
     The aim of this study is to present abstract elements
of a braided  theory  which generalizes standard quantum groups in
a non-trivial and effective way.

     The theory allows  a  possibility  of
completely  pointless  objects  and  includes,  besides   standard
quantum groups, various geometrically interesting structures which
are not quantum groups (Hopf algebras \cite{A}), but which
are more or less similar to them.

     Let us start with a very simple  geometrical  consideration.
Let $G$ be a  quantum  space,  represented  at the formal level
by  a  complex  (unital) *-algebra $\cal{A}$ (consisting  of appropriate
``functions'' on $G$). Geometry of $G$ is encoded in $\cal{A}$.

As an element of the concept of  a  quantum
space, we shall addopt a  point  of  view  according  to  which  $G$
naturally determines  $\cal{A}$,  and  vice  versa.  In  particular,
spaces  with  non-isomorphic but Morita-equivalent ``function''
algebras will be interpreted as  geometrically  different  objects
(although they possess the same basic invariants like  cyclic,  De
Rham and K-homologies \cite{C}).

 If the algebra $\cal{A}$ is noncommutative then the space $G$ can not
be viewed as a structuralized collection of points. However,  this
does not mean that $G$, being a quantum object, possesses no points at
all.

  Actually points of $G$ can be understood as ``maps'' of the form
$\pi\colon\{{\circ}\}\rightarrow G$,
where $\{{\circ}\}$ is the one-point set. Dually, at the level  of
function algebras, every point is determined  by  the  corresponding
``pull back''
$\varphi\colon\cal{A}\rightarrow\Bbb{C}$,
which is a character  (nontrivial  multiplicative  *-functional) on
$\cal{A}$.
Conversely, having Gelfand-Naimark theory of  commutative
$C^*$-algebras in mind (as  well  as  a  possibility  of  recovering
points as characters  in  classical  differential,  or  algebraic,
geometry) it is natural to assume that every character on  $\cal{A}$
determines a point of $G$, in the  way  described.  Of  course,  the
space $G$ may be ``completely quantum''-without points at all.

 Let us now suppose  that  $G$  is  endowed  with  a  quantum  group
structure. By definition, the group structure on $G$ corresponds  to
a Hopf algebra \cite{A} structure  on  $\cal{A}$.  This  structure  is  a
symbiosis of the algebra structure on $\cal{A}$,  and  a  coalgebra
structure specified by the  coproduct
$\phi\colon\cal{A}\rightarrow\cal{A}\otimes
\cal{A}$ and the counit $\e\colon\cal{A}\rightarrow\Bbb{C}$ (we follow the
notation of \cite{W}).
Two structures should  satisfy
certain additional conditions.

     As first, there exists an antipodal map
$\k\colon\cal{A}\rightarrow\cal{A}$ such that
$$m(\k\otimes \id)\phi=m(\id\otimes \k)\phi=1\e,$$
where $m\colon\cal{A}\otimes \cal{A}\rightarrow\cal{A}$
is the multiplication in $\cal{A} $
and $1\in\cal{A}$ is the unit element. Secondly, the map $\phi$ should  be
multiplicative, in the sense that
$$\phi(ab)=\phi(a)\phi(b),\leqno{(*)}$$
for  each  $a,b\in\cal{A}$.  In  the  above  relation,  $\cal{A}\otimes
\cal{A}$ is understood as an algebra, in a natural manner.

 As a consequence of mentioned properties, it turns out  that  the
antipode $\k$ is an anti(co)multiplicative map.
The multiplicativity of the counit is another  important  consequence.
Further, if the group  structure  on  $G$  and  the  *-structure  on
$\cal{A}$ are mutually compatible in the sense that  $\phi*=(*\otimes
*)\phi$ then the composition $*\k$ is involutive  and  the  counit  is
hermitian.
}
 In particular, the space $G$ possesses  at  least  one  point,
corresponding to the counit map (the neutral element). Let  us  denote
by $G_{cl}$   the ``classical part'' of $G$, consisting of all points
of $G$.
The quantum group structure on $G$ induces, in a  natural  manner,  a  group
structure on the set $G_{cl}$   (such that  $G_{cl}$
is  interpretable  as  a subgroup of $G$). As far as
$\cal{A}$  is  noncommutative,  $G_{cl}$    is  a
nontrivial part of $G$. The space $G$ can be imagined as  a  ``disjoint
union'' of two essentially different parts: the classical part
$G_{cl}$ and the purely quantum (pointless) part $G\setminus G_{cl}$.

     In this sense, quantum groups are, in  contrast  to  ordinary
ones, {\it inhomogeneous objects}.

     The mentioned inhomogeneity  is  explicitely  visible  in  the
situations  in  which   ``diffeomorphisms''   of   $G$   appear.   All
diffeomorphisms  must  ``preserve''  the  classical  part  $G_{cl}$.  For
example, in the theory of quantum principal  bundles  over  smooth
manifolds \cite{D} a natural correspondence between quantum  $G$-bundles
and ordinary $G_{cl}$-bundles (over the same manifold) holds.  This  is
because the corresponding right-covariant  ``transition  functions''
are completely determined by their ``restrictions'' on $G_{cl}$  .

     On  the  other  hand,  it  is  natural  to  expect  that   in
noncommutative geometry quantum spaces with a group structure play
a similar role as Lie groups in classical  differential  geometry.
And,  among  other  things,  Lie  groups   provide   examples   of
particularly regular geometrical objects. From this point of view,
the necessity  of  described  classical-quantum  decomposition  of
quantum groups seems strange.

     Such a thinking naturally leads to an idea of generalizing  a
notion of a group structure on a noncommutative space, in order to
include objects of a more elaborate geometrical nature.  This  was
the main motivation for this work.

 Technically speaking,  the  theory  will  be  formulated  in  the
framework of {\it braided algebras}.  Explicitly,  appropriate
linear (braid) operators
$\sigma\colon\cal{A}\otimes   \cal{A}\rightarrow\cal{A}\otimes
\cal{A}$ will enter the game. These operators play the role of  the
standard  transposition  and  induce,  in  a  natural  manner,   a
structure of an associative  algebra  on $\cal{A}\otimes\cal{A}$.
This requires certain compatibility conditions between $\sigma$ and
the algebra $\cal{A}$. We shall also add a  compatibility  condition
between $\sigma$ and the coalgebra  structure (both  compatibility
conditions are trivially fulfilled  in  the  standard  case).  The
generalization will then consist in replacing the axiom $(*)$  by  a
$\sigma$-relativized multiplicativity condition.

     In such a way we obtain ``braided quantum groups''. The attribute
``braided'' will be justified  after  an  appropriate
developement of the formalism, by establishing that $\sigma$  satisfies
the braid equation.

The paper is organized as follows.
The next section is devoted to the definition of braided quantum groups.
In Section 3 the most important interrelations between all  relevant
maps will be investigated. In particular, we shall see that
besides the flip-over  operator  $\sigma$,  another  braid  operator
$\tau\colon\cal{A}\otimes\cal{A}\rightarrow\cal{A}\otimes\cal{A}$   naturally
enters the game. This  operator  is  expressible  via  $\e$, $\phi$  and
$\sigma$.  Two  braid  operators  $\sigma$  and  $\tau$  will   play   a
fundamental role in the whole analysis. In particular, it will  be
shown that $\sigma$ and $\tau$ are mutually compatible in a  ``braided
sense''.

A theory of braided quantum groups presented in
\cite{Maj} can be viewed as a special case of the
theory considered here.

Actually,  braided  quantum  groups introduced in  this
paper reduce to the braided quantum groups of \cite{Maj}
iff two  basic
braid operators coincide. On the  other  hand,  $\sigma=\tau$  is  a
necessary   and    sufficient    condition    for    the
multiplicativity of the counit.
In  particular,  completely  pointless  braided
quantum groups are not includable into the framework of \cite{Maj}.

A large class of examples of completely pointless braided quantum groups
is given by braided Clifford algebras \cite{DD}
associated to involutive braidings.  This includes classical Clifford
and Weyl algebras.

The paper ends with two appendices.
Appendix A is devoted  to  the  main
properties of systems of braid operators, mutually compatible in a
``braided  sense''.  A  motivation  for  that  comes  from   already
mentioned  braided  compatibility  between  $\sigma$  and $\tau$.
In
particular, it will be shown that $\sigma$ and $\tau$ can be naturally
included in a (generally  infinite)  ``braid  system''
expressing concisely all twisting properties. Finally, Appendix B
discusses the question of the multiplicativity of the counit map.

\section{Definition of Braided Quantum Groups}

 Let $\cal{A}$ be a complex associative algebra,  with  the  product
$m\colon\cal{A}\otimes \cal{A}\rightarrow\cal{A}$ and the unit element
$1\in\cal{A}$.  Let  us
assume that $\cal{A}$ is  endowed  with  a  coassociative  coalgebra
structure, specified by the coproduct
$\phi\colon\cal{A}\rightarrow\cal{A}\otimes
\cal{A}$ and the counit $\e\colon\cal{A}\rightarrow\Bbb{C}$.
Finally, let  us  assume  that
bijective linear maps $\k\colon\cal{A}\rightarrow\cal{A}$ and
$\sigma\colon\cal{A}\otimes
\cal{A}\rightarrow\cal{A}\otimes \cal{A}$ are given such that the  following
equalities hold
\begin{align}
\sigma(m\otimes\id)&=(\id\otimes m)(\sigma\otimes\id)(\id\otimes\sigma)
\label{21}\\
\sigma(\id\otimes m)&=(m\otimes\id)(\id\otimes\sigma)(\sigma\otimes \id)
\label{22}\\
\phi m=(m&\otimes m)(\id\otimes\sigma\otimes\id)(\phi\otimes\phi)\label{23}\\
(\sigma\otimes\id^2)(\id\otimes\phi\otimes\id)(&\sigma^{-1}\otimes\id)
(\id\otimes\phi)\label{24}\\
&=(\id^2\otimes\sigma)(\id\otimes\phi\otimes\id)(\id\otimes
\sigma^{-1})(\phi\otimes\id),\notag
\end{align}
together with the antipode axiom
\begin{equation}
1\e=m(\id\otimes\k)\phi=m(\k\otimes\id)\phi.\label{25}
\end{equation}

\begin{defn}\label{def:21} Every pair
$G=\bigl(\cal{A},\{\phi,\e,\k,\sigma\}\bigr)$
satisfying the above requirements is called {\it a braided quantum group}.
\end{defn}

 The map $\sigma$ is interpretable as the  ``twisting  operator''.  In
the standard theory, $\sigma$ reduces to the ordinary transposition.
Identities \eqref{21}--\eqref{24} express  mutual
compatibility  between
maps $\phi$, $m$ and $\sigma$.  It  is  important  to  mention  that  the
asymmetry between \eqref{21}--\eqref{22} and \eqref{24}
implies that the theory is not
``selfdual''.  However if  we  replace \eqref{24}   with   ``dual''
counterparts of \eqref{21}--\eqref{22}
then the theory reduces to braided quantum
groups of \cite{Maj} (and in particular becomes selfdual).

 The space $\cal{A}\otimes \cal{A}$ is  an  $\cal{A}$-bimodule,  in  a
natural manner. With the help of $\sigma$, a natural product can  be
defined on $\cal{A}\otimes \cal{A}$, by requiring
\begin{equation}
 (a\otimes b)(c\otimes d)=a\sigma(b\otimes c)d.\label{26}
\end{equation}
Identities \eqref{21}--\eqref{22} ensure that this defines
an  associative  algebra
structure on $\cal{A}\otimes \cal{A}$, such that $1\otimes 1$  is  the
unit element. In particular,
\begin{equation}\label{27}
\sigma\bigl(1\otimes (\,)\bigr)=(\,)\otimes 1 \qquad\quad
\sigma\bigl((\,)\otimes 1\bigr)=1\otimes (\,).
\end{equation}
In the following, it will be assumed that  $\cal{A}\otimes  \cal{A}$
is endowed with this algebra structure.  Equality \eqref{23}  then
says that $\phi$ is multiplicative.

Identity \eqref{24}  expresses  the  coassociativity  of  the  map
$(\id\otimes \sigma^{-1}  \otimes \id)(\phi\otimes  \phi)$.
The  ``inverse'' identity
\begin{multline}\label{28}
(\id^2\otimes\sigma)(\id\otimes\phi\otimes\id)(\sigma\otimes\id)
(\id\otimes\phi)\\
=(\sigma\otimes\id^2)(\id\otimes\phi\otimes\id)
(\id\otimes
\sigma)(\phi\otimes\id)
\end{multline}
holds, too.
It expresses  the  coassociativity  of $(\id\otimes  \sigma\otimes
\id)(\phi\otimes \phi)$.

\section{Elementary Algebraic Properties}

Let $G=\bigl(\cal{A},\{\phi,\e,\k,\sigma\}\bigr)$ be a
braided quantum group. As in the standard theory, the
antipode is uniquely determined by
\eqref{25}. The flip-over operator $\sigma$ is expressible through
$\phi,m$ and $\k$ in the following way
\begin{equation}\label{29}
\sigma=(m\otimes m)(\k\otimes\phi m\otimes\k)(\phi\otimes \phi),
\end{equation}
as directly follows from \eqref{23} and \eqref{25}.

     It is easy to see that
\begin{equation}
\phi(1)=1\otimes 1.\label{210}
\end{equation}
Indeed, $\phi(1)$    is    the    unity in the    subalgebra
$\phi(\cal{A})\subseteq\cal{A}\otimes \cal{A}$, as  follows
from \eqref{23}. On the other hand, $\cal{A}\otimes \cal{A}$ is generated by
$\phi(\cal{A})$, as a left
(right) $\cal{A}$-module. Hence, $\phi(1)$ is the unity of $\cal{A}\otimes
\cal{A}$. From \eqref{210} we obtain
\begin{align}
                            \e(1)&=1\label{211} \\
                            \k(1)&=1.\label{212}
\end{align}
 In   further   computations   the   result   of   an $(n-1)$-fold
comultiplication of an  element  $a\in\cal{A}$  will  be  symbolically
denoted by $a^{(1)}\otimes \dots\otimes a^{(n)}$.
Clearly, this  element  of $\cal{A}$   is  independent  of
ways  in  which  the  corresponding
comultiplications are performed.
\begin{lem}\label{lem:21}
The following identities hold
\begin{align}
 (\e\otimes\id)&=(\id\otimes
\e m)(\sigma\otimes  \id)(\id\otimes \phi) \label{213} \\
(\id\otimes \e)&=(\e m\otimes\id)(\id\otimes\sigma)(\phi\otimes\id).\label{214}
\end{align}
\end{lem}
\begin{pf} According to \eqref{23},
$$ab^{(1)}\otimes b^{(2)}=(\e\otimes\id)(m\otimes m)(\id\otimes
\sigma\otimes \id)\bigl(
a^{(1)}\otimes a^{(2)}   \otimes b^{(1)}   \otimes b^{(2)}\bigr)\otimes
b^{(3)},$$
for each $a,b\in\cal{A}$. Acting by $m(\id\otimes \k)$ on  this  equality,
and using \eqref{25} we obtain
$$ a\e(b)=(\e\otimes \id)\bigl(a^{(1)}\sigma(a^{(2)}\otimes b)\bigr).$$
Similarly, acting by $m(\k\otimes \id)$ on the identity
$$
a^{(1)}\otimes a^{(2)}b=a^{(1)}   \otimes (\id\otimes \e)
(m\otimes m)(\id\otimes \sigma\otimes \id)(a^{(2)}\otimes
a^{(3)}\otimes b^{(1)}   \otimes b^{(2)}   ) $$
we obtain
$$\e(a)b=(\id\otimes \e)\bigl(\sigma(a\otimes b^{(1)})b^{(2)}\bigr).\qed$$
\renewcommand{\qed}{}
\end{pf}

 A ``secondary'' flip-over operator $\tau$ will be now introduced in the
game. From \eqref{24} we obtain
\begin{equation}
(\id^2 \otimes \e)(\id\otimes \sigma^{-1}  )(\phi\otimes \id)=
(\e\otimes \id^2 )(\sigma^{-1}\otimes \id)(\id\otimes \phi). \label{215}
\end{equation}
     Let $\tau\colon\cal{A}\otimes \cal{A}\rightarrow\cal{A}\otimes \cal{A}$
be a linear map defined by
\begin{equation}
\tau=(\id^2 \otimes \e)(\id\otimes\sigma^{-1})(\phi\otimes
\id)\sigma=(\e\otimes \id^2)(
\sigma^{-1}\otimes \id)(\id\otimes\phi)\sigma. \label{216}
\end{equation}
\begin{lem}\label{lem:22}
The map $\tau$ is bijective and
\begin{equation}
\tau^{-1}\sigma=(\id^2 \otimes \e)(\id\otimes \sigma)(\phi\otimes
\id)=(\e\otimes\id^2)
(\sigma\otimes\id)(\id\otimes \phi).\label{217}
\end{equation}
\end{lem}
\begin{pf}
The second equality in \eqref{217} follows from \eqref{28}.
Let $\tau'\sigma$ be the map given by the second term in \eqref{217}.  A
direct computation gives
\begin{equation*}
\begin{split}
{}&\tau\tau'\sigma=(\e\otimes \id^2 \otimes \e)(\sigma^{-1}\otimes \id^2)
(\id\otimes \phi
\otimes\id)(\sigma\otimes\id)(\id\otimes\sigma)(\phi\otimes\id)\\
=&(\e\otimes    \id^2 \otimes     \e)(\sigma^{-1}  \otimes\id^2 )
(\id\otimes \phi\otimes \id)(\sigma\otimes \e\otimes \id)
(\id\otimes \phi\otimes \id)(\id\otimes \sigma)(\phi\otimes\id)\\
=&(\e\otimes\id^2 \otimes   \e\otimes   \e)(\sigma^{-1}  \otimes   \id
\otimes \sigma)(\id\otimes\phi\otimes\id^2)(\id\otimes\phi\otimes
\id)(\sigma\otimes \id)(\id\otimes\phi)\\
=&(\e\otimes \id^2 \otimes  \e\otimes  \e)(\id^3 \otimes\sigma)(\id^2
\otimes \phi\otimes\id)(\id^2 \otimes\sigma^{-1})(\id\otimes\phi\otimes
\id)(\id\otimes \sigma)(\phi\otimes \id)\\
=&(\id^2 \otimes \e\otimes \e)(\id^2 \otimes \sigma)(\id\otimes  \phi\otimes
\id)(\id\otimes \sigma^{-1}  )(\phi\otimes \id)\sigma\\
=&(\id^2 \otimes \e\otimes \e)(\sigma\otimes\id^2 )(\id\otimes  \phi\otimes
\id)(\sigma^{-1}  \otimes \id)(\id\otimes \phi)\sigma=\sigma.
\end{split}
\end{equation*}

 Similarly,  interchanging  $\sigma$  and  $\sigma^{-1}$    in  the   above
computations we conclude that $\tau'$ is a left inverse for
$\tau$.  Hence,  $\tau$  is bijective and $\tau^{-1}=\tau'$.
\end{pf}

 Let us write down some important  algebraic  relations  including
the map $\tau$. As first, let us observe that
\begin{gather}
(\e\otimes \id)\tau=\id\otimes \e \quad\qquad
(\id\otimes \e)\tau=\e\otimes \id  \label{218}\\
\tau\bigl(1\otimes (\,)\bigr)=(\,)\otimes 1\quad\qquad
\tau\bigl((\,)\otimes 1\bigr)=1\otimes (\,). \label{219}
\end{gather}
This is a direct consequence of the definition of $\tau$, and  property
\eqref{27}. Further, coassociativity of $\phi$ and relations
\eqref{216}--\eqref{217} imply
\begin{align}
(\phi\otimes\id)\tau^{-1}\sigma&
=(\id\otimes\tau^{-1}\sigma)(\phi\otimes\id)\label{220}\\
(\id\otimes\phi)\tau^{-1}\sigma&=(\tau^{-1}\sigma\otimes\id)(\id\otimes\phi)
\label{221}\\
(\phi\otimes\id)\tau\sigma^{-1}&=(\id\otimes\tau\sigma^{-1})(\phi\otimes\id)
\label{222}\\
(\id\otimes \phi)\tau\sigma^{-1}&=(\tau\sigma^{-1}\otimes \id)
(\id\otimes \phi).\label{223}
\end{align}
In  other  words,   maps   $\sigma\tau^{-1}$     and
$\sigma^{-1}\tau$   are
automorphisms of  the  $\cal{A}$-bicomodule  $\cal{A}\otimes  \cal{A}$
(with the left and the right $\cal{A}$-comodule structures given  by
$\phi\otimes \id$  and  $\id\otimes  \phi$  respectively).  In  general,
$\sigma\tau^{-1}$    and  $\sigma^{-1}\tau$  do  not  commute.
However, the following commutation relations hold
\begin{align}
(\sigma\tau^{-1}\otimes \id)(\id\otimes \sigma\tau^{-1})&
=(\id\otimes\sigma\tau^{-1})(\sigma\tau^{-1}\otimes \id)\label{224}\\
(\sigma\tau^{-1}\otimes \id)(\id\otimes \sigma^{-1}\tau)&
=(\id\otimes \sigma^{-1}\tau)(\sigma\tau^{-1}\otimes \id) \label{225} \\
(\sigma^{-1}\tau\otimes \id)(\id\otimes \sigma\tau^{-1})&
=(\id\otimes \sigma\tau^{-1})(\sigma^{-1}\tau\otimes \id)\label{226}\\
(\sigma^{-1}\tau\otimes \id)(\id\otimes \sigma^{-1}\tau)&
=(\id\otimes \sigma^{-1}\tau)(\sigma^{-1}\tau\otimes \id).\label{227}
\end{align}

The above equalities follow from \eqref{220}--\eqref{223} and
\eqref{216}--\eqref{217}. As a direct consequence of Lemma~\ref{lem:21} and
\eqref{217} we find
\begin{equation}
\e m=(\e\otimes \e)\sigma^{-1}\tau. \label{228}
\end{equation}
This generalizes the standard multiplicativity law for the counit.

Identities \eqref{24} and \eqref{28} can be rewritten in a simpler
``pentagonal form'', including the  operator
$\tau$  and  explicitly expressing  twisting properties of the coproduct map.
\begin{pro}\label{pro:23}
The following identities hold
\begin{align}
(\phi\otimes\id)\sigma&=(\id\otimes\tau)(\sigma\otimes\id)(\id\otimes\phi)
\label{229}\\
(\id\otimes\phi)\sigma&=(\tau\otimes\id)(\id\otimes\sigma)(\phi\otimes\id)
\label{230}\\
(\phi\otimes\id)\sigma&=(\id\otimes\sigma)(\tau\otimes\id)(\id\otimes\phi)
\label{231}\\
(\id\otimes\phi)\sigma&=(\sigma\otimes\id)(\id\otimes\tau)(\phi\otimes\id).
\label{232}
\end{align}
\end{pro}
\begin{pf}
Using \eqref{24} and \eqref{217} we obtain
\begin{multline*}
(\e\otimes \id^3 )(\sigma\otimes \id^2 )(\id\otimes \phi\otimes \id)
(\sigma^{-1}\otimes\id)(\id\otimes\phi)=(\tau^{-1}\otimes\id)(\id\otimes\phi)\\
=(\e\otimes \id\otimes \sigma)(\id\otimes \phi\otimes \id)
(\id\otimes \sigma^{-1})(\phi\otimes \id)=(\id\otimes \sigma)(\phi\otimes \id)
\sigma^{-1}.
\end{multline*}
Similarly,
\begin{multline*}
(\id^3 \otimes    \e)(\id^2 \otimes\sigma)(\id\otimes\phi\otimes
\id)(\id\otimes \sigma^{-1})(\phi\otimes \id)=(\id\otimes \tau^{-1})
(\phi\otimes \id)\\
=(\sigma\otimes\id\otimes\e)(\id\otimes\phi\otimes\id)(\sigma^{-1}\otimes
\id)(\id\otimes\phi)=(\sigma\otimes\id)(\id\otimes \phi)\sigma^{-1}.
\end{multline*}
     Hence, \eqref{229}--\eqref{230}  hold.  Starting  from
equalities \eqref{28}  and \eqref{216}  and  applying  the  same
computation we obtain \eqref{231}--\eqref{232}.
\end{pf}

 In the next proposition ``pentagonal''  twisting  relations  including
only $\tau$ are collected.
\begin{pro}\label{pro:24} We have
\begin{align}
(\phi\otimes\id)\tau&=(\id\otimes\tau)(\tau\otimes\id)(\id\otimes\phi)
\label{233}\\
(\id\otimes\phi)\tau&=(\tau\otimes\id)(\id\otimes\tau)(\phi\otimes\id)
\label{234}\\
\tau(m\otimes\id)&=(\id\otimes\tau)(\tau\otimes\id)(\id\otimes m)
\label{235}\\
\tau(\id\otimes m)&=(\tau\otimes\id)(\id\otimes\tau)(m\otimes \id).
\label{236}
\end{align}
\end{pro}
\begin{pf} Direct transformations give
$$(\id\otimes\tau)(\tau\otimes\id)(\id\otimes\phi)=(\id\otimes
\tau\sigma^{-1})(\phi\otimes \id)\sigma=(\phi\otimes \id)\tau.$$\
Similarly,
$$(\tau\otimes\id)(\id\otimes\tau)(\phi\otimes
\id)=(\tau\sigma^{-1}\otimes\id)(\id\otimes\phi)\sigma=(\id\otimes
\phi)\tau.$$
     Applying \eqref{216}, \eqref{231} and \eqref{21} we obtain
\begin{equation*}
\begin{split}
(\id\otimes m)&(\tau\otimes \id)(\id\otimes \tau)=(\id\otimes  m\otimes
\e)(\tau\otimes\sigma^{-1})(\id\otimes\phi\otimes\id)(\id\otimes
\sigma)\\
&=(\id\otimes m\otimes\e)(\id^2\otimes\sigma^{-1})(\id\otimes
\sigma^{-1}\otimes\id)(\phi\otimes  \id^2 )(\sigma\otimes  \id)(\id\otimes
\sigma)\\
&=(\id^2 \otimes \e)(\id\otimes\sigma^{-1})(\phi\otimes  m)(\sigma\otimes
\id)(\id\otimes\sigma)\\
&=(\id^2\otimes\e)(\id\otimes
\sigma^{-1})(\phi\otimes \id)\sigma(m\otimes \id)=
\tau(m\otimes \id).
\end{split}
\end{equation*}
Similarly,
\begin{equation*}
\begin{split}
(m\otimes \id)&(\id\otimes \tau)(\tau\otimes\id)=(\e\otimes  m\otimes
\id)(\sigma^{-1}\otimes \tau)(\id\otimes \phi\otimes  \id)(\sigma\otimes
\id)\\
&=(\e\otimes m\otimes\id)(\sigma^{-1}\otimes\id^2)(\id\otimes
\sigma^{-1}\otimes\id)(\id^2 \otimes\phi)(\id\otimes
\sigma)(\sigma\otimes \id)\\
&=(\e\otimes   \id^2 )(\sigma^{-1}\otimes\id)(m\otimes\phi)(\id\otimes
\sigma)(\sigma\otimes \id)=\tau(\id\otimes m).\qed
\end{split}
\end{equation*}
\renewcommand{\qed}{}
\end{pf}
     We pass to the study of  algebraic  relations  including  the
antipode map. In the standard theory, the antipode is an
anti(co)-multiplicative map. The next proposition  gives  a  braided
counterpart of this property.
\begin{pro}\label{pro:25}
We have
\begin{align}
\phi\k&=\sigma(\k\otimes \k)\phi  \label{237}\\
\k m&=m(\k\otimes \k)\tau\sigma^{-1}\tau\sigma^{-1}\tau. \label{238}
\end{align}
\end{pro}
\begin{pf}
Let us start from the identity
$$ \k(a^{(1)})a^{(2)}\otimes a^{(3)}=1\otimes a. $$
Acting by $\phi\otimes \phi$ on both sides, and using \eqref{23} and
\eqref{210} we obtain
$$\bigl(\phi\k(a^{(1)})\bigr)\bigl(a^{(2)}\otimes a^{(3)}\bigr)
\otimes a^{(4)}\otimes a^{(5)}=1\otimes 1\otimes a^{(1)}\otimes a^{(2)}.$$
After the action of $(\id\otimes m\otimes \id)(\id^2\otimes \k\otimes \id)$
on both sides the above equality becomes
$$ \bigl(\phi\k(a^{(1)})\bigr)\bigl(a^{(2)}\otimes
1\bigr)\otimes a^{(3)}=1\otimes
\k(a^{(1)})\otimes a^{(2)}.$$
Hence
$$
\bigl(\phi\k(a^{(1)})\bigr)\bigl(a^{(2)}\k(a^{(3)})\otimes 1\bigr)
=\bigl(1\otimes \k(a^{(1)})\bigr)\bigl(\k(a^{(2)})\otimes 1\bigr).
$$
Applying \eqref{26}--\eqref{27} we obtain
$$\phi\k(a)=\sigma\bigl(\k(a^{(1)})\otimes \k(a^{(2)})\bigr). $$
This proves \eqref{237}.  To  prove  \eqref{238},  let  us  start  from
$m(\k\otimes m)(\phi\otimes \id)=\e\otimes \id$, act by  it  on
$m\otimes m$, and apply \eqref{23} and \eqref{228}. We find
$$
m(\k\otimes m)(m\otimes m\otimes m)(\id\otimes \sigma\otimes \id^3 )(\phi
\otimes\phi\otimes \id^2 )=(\e\otimes \e)\sigma^{-1}\tau\otimes m. $$
Acting by this  equality  on  $(\id^2 \otimes  \k\otimes  \id)(\id\otimes
\phi\otimes \id)$ and simplifying the expression we find
$$
m(\k m\otimes m)(\id\otimes     \sigma\otimes     \id)(\phi\otimes
\id^2)=(\e\otimes m)(\id\otimes \k\otimes \id)(\sigma^{-1}\tau\otimes \id).$$
Acting by this on $(\id^2\otimes
\k)(\id\otimes \sigma)(\phi\otimes \id)$ we obtain
\begin{multline*}
m(\k m\otimes m)(\id\otimes\sigma\otimes\k)(\id^2 \otimes
\sigma)(\id\otimes \phi\otimes \id)(\phi\otimes \id)\\
=(\e\otimes m)(\id\otimes \k\otimes \k)(\sigma^{-1}\tau\otimes
\id)(\id\otimes \sigma)(\phi\otimes \id).
\end{multline*}
After simple twisting transformations the left-hand  side  of  the
above equality becomes
\begin{multline*}
m(\k m\otimes m)(\id^3 \otimes \k)(\id^2 \otimes \phi)(\id\otimes
\sigma\tau^{-1}  \sigma
)(\phi\otimes \id)\\
=(\k m\otimes \e)(\id\otimes \sigma\tau^{-1}\sigma)(\phi\otimes \id)
=\k m\tau^{-1}\sigma\tau^{-1}\sigma.
\end{multline*}
The right-hand side of the mentioned equality reduces to
$$m(\e\otimes\k\otimes\k)(\sigma^{-1}\otimes\id)(\id\otimes
\phi)\sigma=m(\k\otimes\k)\tau.$$
Consequently, \eqref{238} holds.
\end{pf}

     Twisting properties of the antipode will be now analyzed.  As
first, a technical lemma
\begin{lem}\label{lem:26} We have
\begin{gather}
\bigl[\sigma(\k\otimes \id)\tau^{-1}\sigma\tau^{-1}(a\otimes b^{(1)})\bigr]
b^{(2)}=a\otimes 1\e(b)\\
a^{(1)}\bigl[\sigma(\id\otimes \k)\tau^{-1}\sigma\tau^{-1}(a^{(2)}\otimes b)
\bigr]=\e(a)1\otimes b,
\end{gather}
for each $a,b\in\cal{A}$.
\end{lem}
\begin{pf}
We compute
\begin{equation*}
\begin{split}
&(\id\otimes m)(\sigma\otimes\id)(\k\otimes
\id^2)(\tau^{-1}\sigma\tau^{-1}\otimes\id)(\id\otimes \phi)\\
=&(\id\otimes   m)(\sigma\otimes\id)(\k\otimes   \id^2 )(\tau^{-1}\otimes
\id)(\id\otimes\phi)\sigma\tau^{-1}\\
=&(\id\otimes  m)(\sigma\otimes   \id)(\k\otimes   \sigma)(\phi\otimes
\id)\tau^{-1}=\sigma(m\otimes \id)(\k\otimes \id^2)(\phi\otimes \id)\tau^{-1}\\
=&\sigma(1\e\otimes \id)\tau^{-1}=\id\otimes 1\e.
\end{split}
\end{equation*}
Similarly,
\begin{equation*}
\begin{split}
&(m\otimes \id)(\id\otimes \sigma)(\id^2 \otimes\k)(\id\otimes \tau^{-1}
\sigma\tau^{-1})(\phi\otimes \id)\\
=&(m\otimes    \id)(\id\otimes    \sigma)(\id^{2}\otimes\k)(\id\otimes
\tau^{-1})(\phi\otimes\id)\sigma\tau^{-1}\\
=&\sigma(\id\otimes m)(\id \otimes\k)(\id\otimes
\phi)\tau^{-1}=\sigma(\id\otimes 1\e)\tau^{-1}=1\e\otimes\id.\qed
\end{split}
\end{equation*}
\renewcommand{\qed}{}
\end{pf}
\begin{pro}\label{pro:27}
The following identities hold
\begin{gather}
\sigma(\k\otimes \id)=(\id\otimes \k)\tau\sigma^{-1}\tau\label{242}\\
\tau(\id\otimes \k)=(\k\otimes \id)\tau\label{243}\\
\tau(\k\otimes \id)=(\id\otimes \k)\tau\label{244}\\
\sigma(\id\otimes\k)=(\k\otimes \id)\tau\sigma^{-1}\tau.\label{241}
\end{gather}
\end{pro}
\begin{pf}
Applying Lemma~\ref{lem:26} and property \eqref{25} we obtain
$$\sigma(\k\otimes  \id)\tau^{-1}\sigma\tau^{-1}(a\otimes b)=
\bigl[\sigma(\k\otimes \id)\tau^{-1}\sigma\tau^{-1}(a\otimes b^{(1)})\bigr]
b^{(2)}\k(b^{(3)})=a\otimes \k(b). $$
Similarly,
$$(\id\otimes \k)\tau^{-1}\sigma\tau^{-1}(a\otimes b)
=\k(a^{(1)})a^{(2)}\bigl[\sigma(\k\otimes    \id)\tau^{-1}\sigma\tau^{-1}
(a^{(3)}   \otimes b)\bigr]=\k(a)\otimes b.$$
This shows \eqref{242} and \eqref{241}.
Using properties \eqref{216}, \eqref{242}, \eqref{241},
\eqref{222}--\eqref{223}  and  \eqref{233}--\eqref{234}
we obtain
\begin{equation*}
\begin{split}
\tau(\id\otimes \k)&=(\e\otimes \id^2)(\sigma^{-1}\otimes \id)
(\id\otimes \phi)\sigma(\id\otimes \k)\\
&=(\e\otimes \k\otimes \id)(\tau^{-1}\sigma\tau^{-1}\otimes \id)
(\id\otimes \phi)\tau\sigma^{-1}\tau\\
&=(\e\otimes\k\otimes \id)(\tau^{-1}\otimes \id)(\id\otimes \phi)\tau=(\k
\otimes\id)\tau.
\end{split}
\end{equation*}
Similarly,
\begin{equation*}
\begin{split}
\tau(\k\otimes \id)&=(\id^2 \otimes \e)(\id\otimes \sigma^{-1})
(\phi\otimes \id)\sigma(\k\otimes \id)\\
&=(\id\otimes\k\otimes\e)(\id\otimes \tau^{-1}\sigma\tau^{-1})
(\phi\otimes \id)\tau\sigma^{-1}\tau\\
&=(\id\otimes\k\otimes \e)(\id\otimes \tau^{-1})(\phi\otimes \id)\tau=
(\id\otimes\k)\tau. \qed
\end{split}
\end{equation*}
\renewcommand{\qed}{}
\end{pf}

     As a direct consequence of the previous proposition we find
\begin{align}
(\k\otimes\k)\tau&=\tau(\k\otimes \k)\label{245}\\
(\k\otimes\k)\sigma&=\sigma(\k\otimes\k).\label{246}
\end{align}
For the end of this section, we shall prove that $\sigma$ and  $\tau$
satisfy a system of braid equations.
\begin{pro}\label{pro:28}
The following identities hold
\begin{align}
(\sigma\otimes \id)(\id\otimes \sigma)(\sigma\otimes \id)&=(\id\otimes
\sigma)(\sigma\otimes \id)(\id\otimes \sigma)\label{247}\\
(\tau\otimes \id)(\id\otimes \sigma)(\sigma\otimes \id)&=(\id\otimes \sigma)(
\sigma\otimes \id)(\id\otimes \tau)\label{248}\\
(\sigma\otimes \id)(\id\otimes \tau)(\sigma\otimes \id)&=(\id\otimes \sigma
)(\tau\otimes \id)(\id\otimes \sigma)\label{249} \\
(\sigma\otimes \id)(\id\otimes \sigma)(\tau\otimes \id)&=(\id\otimes \tau)(
\sigma\otimes \id)(\id\otimes \sigma)\label{250}\\
(\tau\otimes \id)(\id\otimes \tau)(\sigma\otimes\id)
&=(\id\otimes \sigma)(\tau\otimes \id)(\id\otimes \tau)\label{251}\\
(\tau\otimes \id)(\id\otimes \sigma)(\tau\otimes
\id)&=(\id\otimes \tau)(\sigma\otimes \id)(\id\otimes \tau)\label{252}\\
(\sigma\otimes \id)(\id\otimes \tau)(\tau\otimes \id)&=(\id\otimes \tau)(
\tau\otimes \id)(\id\otimes \sigma)\label{253}\\
(\tau\otimes \id)(\id\otimes \tau)(\tau\otimes \id)&=(\id\otimes
\tau)(\tau\otimes \id)(\id\otimes\tau).\label{254}
\end{align}
\end{pro}
\begin{pf} We shall first prove \eqref{248}--\eqref{251}  and  \eqref{253},
secondly
\eqref{254}, thirdly \eqref{252} and finally \eqref{247}. A direct
computation gives
\begin{multline*}
(\tau\otimes \id)(\id\otimes \sigma)(\sigma\otimes  \id)=(\tau\otimes
\id)(\id\otimes \sigma)(m\otimes m\otimes \id)\\
\hfill(\k\otimes\phi m\otimes
\k\otimes \id)(\phi\otimes \phi\otimes \id)\\
=(\id\otimes  m\otimes  m)(\tau\otimes  \id^3)(\id\otimes\tau\otimes
\id^2)(\id^2\otimes\sigma\otimes\id)(\id^3\otimes\sigma)\hfill\\
\hfill(\k\otimes
\phi m\otimes \k\otimes\id)(\phi\otimes \phi\otimes \id)\\
=(\id\otimes    m\otimes    m)(\id\otimes    \k\otimes    \phi\otimes
\k)(\tau\otimes \id^2)(\id\otimes \sigma\otimes \id)(\id\otimes m\otimes
\tau\sigma^{-1}\tau)(\phi\otimes \phi\otimes \id)\hfill\\
=(\id\otimes m\otimes    m)(\id\otimes\k\otimes\phi m\otimes
\k)(\tau\otimes \id^3 )(\id\otimes\sigma\otimes\id^2)(\phi\otimes
\id\otimes \phi)(\id\otimes \tau)\hfill\\
=(\id\otimes   m\otimes    m)(\id\otimes    \k\otimes    \phi m\otimes
\k)(\id\otimes   \phi\otimes    \phi)(\sigma\otimes    \id)(\id\otimes
\tau)\hfill\\
=(\id\otimes \sigma)(\sigma\otimes \id)(\id\otimes \tau).
\end{multline*}
Similarly,
\begin{multline*}
(\id\otimes \tau)(\sigma\otimes \id)(\id\otimes \sigma)=(\id\otimes \tau)(
\sigma\otimes   \id)(\id\otimes   m\otimes   m)\\
\hfill(\id\otimes   \k\otimes
\phi m\otimes \k)(\id\otimes \phi\otimes \phi)\\
=(m\otimes m\otimes \id)(\id^3 \otimes \tau)(\id^2 \otimes \tau\otimes \id)
(\id\otimes\sigma\otimes \id^2)(\sigma\otimes \id^3)\hfill\\
\hfill(\id\otimes \k
\otimes \phi m\otimes\k)(\id\otimes \phi\otimes \phi)\\
=(m\otimes m\otimes \id)(\k\otimes\phi\otimes\k\otimes
\id)(\id^2 \otimes\tau)(\id\otimes\sigma\otimes
\id)(\tau\sigma^{-1}\tau\otimes  m\otimes  \id)(\id\otimes   \phi\otimes
\phi)\hfill\\
=(m\otimes m\otimes \id)(\k\otimes \phi m\otimes \k\otimes \id)
(\id^3 \otimes \tau)(\id^2 \otimes\sigma\otimes\id)(\phi\otimes\id\otimes
\phi)(\tau\otimes \id)\hfill\\
=(m\otimes m\otimes\id)(\k\otimes \phi m\otimes\k\otimes
\id)(\phi\otimes  \phi\otimes   \id)(\id\otimes   \sigma)(\tau\otimes
\id)\hfill\\
=(\sigma\otimes \id)(\id\otimes \sigma)(\tau\otimes \id).
\end{multline*}
Essentially  the  same  transformations  lead  to  identities
\eqref{249}, \eqref{251} and \eqref{253}. Let us prove \eqref{254}.
We have
\begin{equation*}
\begin{split}
&(\id\otimes   \tau)(\tau\otimes   \id)(\id\otimes    \tau)=(\id\otimes
\tau\otimes\e)(\tau\otimes     \sigma^{-1})(\id\otimes\phi\otimes
\id)(\id\otimes \sigma)\\
=&(\id^2 \otimes \e\otimes \id)(\id^2 \otimes  \tau)(\id\otimes  \tau\otimes
\id)(\tau\otimes  \sigma^{-1})(\id\otimes   \phi\otimes   \id)(\id\otimes
\sigma)\\
=&(\id^2 \otimes \e\otimes \id)(\id\otimes \sigma^{-1}\otimes \id)
(\id^2 \otimes \tau)(\id\otimes \tau\otimes \id)(\tau\otimes \id^2)
(\id\otimes \phi\otimes \id)(\id\otimes\sigma)\\
=&(\id^2 \otimes\e\otimes\id)(\id\otimes \sigma^{-1}\otimes
\id)(\phi\otimes \tau)(\tau\otimes \id)(\id\otimes \sigma)\\
=&(\tau\sigma^{-1}\otimes \id)(\id\otimes \tau)(\tau\otimes \id)(\id
\otimes\sigma)
=(\tau\otimes \id)(\id\otimes \tau)(\tau\otimes \id).
\end{split}
\end{equation*}
Identities \eqref{225}, \eqref{248}, \eqref{251} and \eqref{254} imply
\begin{equation*}
\begin{split}
(\id\otimes\tau)&(\sigma\otimes\id)(\id\otimes\tau)=(\id\otimes
\tau)(\sigma\tau^{-1}\otimes \id)(\tau\otimes \id)(\id\otimes\tau)\\
=&(\id\otimes\sigma)(\sigma\tau^{-1}\otimes\id)(\id\otimes
\sigma^{-1}\tau)(\tau\otimes \id)(\id\otimes \tau)\\
=&(\id\otimes\sigma)(\sigma\tau^{-1}\otimes\id)(\id\otimes
\sigma^{-1})(\tau\otimes \id)(\id\otimes \tau)(\tau\otimes \id)\\
=&(\id\otimes\sigma)(\sigma\otimes\id)(\id\otimes
\tau)(\sigma^{-1}\tau\otimes\id)=(\tau\otimes\id)(\id\otimes
\sigma)(\tau\otimes \id).
\end{split}
\end{equation*}
     Finally, \eqref{224}, \eqref{248}, \eqref{250} and \eqref{252}
imply
\begin{equation*}
\begin{split}
(\id\otimes \sigma)&(\sigma\otimes \id)(\id\otimes\sigma)=(\id\otimes
\sigma\tau^{-1})(\sigma\otimes \id)(\id\otimes \sigma)(\tau\otimes \id)\\
&=(\id\otimes \sigma\tau^{-1})(\sigma\tau^{-1}\otimes \id)(\id\otimes
\tau)(\sigma\otimes \id)(\id\otimes \tau)\\
&=(\sigma\tau^{-1}\otimes \id)(\id\otimes \sigma)(\sigma\otimes \id)
(\id\otimes\tau)
=(\sigma\otimes \id)(\id\otimes \sigma)(\sigma\otimes \id).\qed
\end{split}
\end{equation*}
\renewcommand{\qed}{}
\end{pf}
\appendix
\section{Braid Systems}
Let us consider a complex associative algebra $\cal{A}$ with the unit element
$1\in\cal{A}$ and the product
$m\colon\cal{A}\otimes\cal{A}\rightarrow\cal{A}$.
\begin{defn}\label{def:A1} A {\it braid system} over $\cal{A}$ is a collection
$\cal{F}$
of bijective linear maps acting in $\cal{A}\otimes \cal{A}$ and satisfying
\begin{align}
(\alpha\otimes\id)(\id\otimes\beta)(\gamma\otimes\id)&=
(\id\otimes\gamma)(\beta\otimes\id)(\id\otimes\alpha)\label{A1}\\
\alpha(\id\otimes m)&=(m\otimes\id)(\id\otimes\alpha)
(\alpha\otimes\id)\label{A2}\\
\alpha(m\otimes\id)&=(\id\otimes m)(\alpha\otimes\id)
(\id\otimes\alpha)\label{A3}
\end{align}
for each $\alpha,\beta,\gamma\in\cal{F}$.
\end{defn}
\begin{defn}\label{def:A2}
A braid system $\cal{F}$ is called {\it complete} iff it is closed under the
operation $(\alpha,\beta,\gamma)\mapsto\alpha\beta^{-1}\gamma$.
\end{defn}
Let $\cal{F}$ be a braid system over $\cal{A}$. Then
$$
\alpha\bigl(1\otimes (\,)\bigr)=(\,)\otimes 1\quad\qquad
\alpha\bigl((\,)\otimes 1\bigr)=(\,)\otimes 1
$$
for each $\alpha\in\cal{F}$, as follows from \eqref{A2}--\eqref{A3}.
Further, every $\alpha\in\cal{F}$ naturally determines an associative algebra
structure on $\cal{A}\otimes\cal{A}$, with the unit element $1\otimes 1$.
The corresponding product is given by
$(m\otimes m)(\id\otimes\alpha\otimes\id)$.

We are going to prove that there exists the {\it minimal} complete
braid system $\cal{F}^*$ which extends $\cal{F}$.
Starting from the system $\cal{F}$ we can inductively
construct an increasing chain of
braid systems $\cal{F}_n$, where $n\geq 0$ and
$\cal{F}_0=\cal{F}$, while $\cal{F}_{n+1}$ is consisting of maps of the form
$\delta=\alpha\beta^{-1}\gamma$, where $\alpha,\beta,\gamma\in\cal{F}_n$.
The fact
that all $\cal{F}_n$ are braid systems easily follows by
induction, applying the definition of braid systems and the identity
\begin{equation}
(\alpha\beta^{-1}\otimes\id)(\id\otimes\gamma\delta^{-1})
=(\id\otimes\gamma\delta^{-1})(\alpha\beta^{-1}\otimes\id)\label{A4}
\end{equation}
(which holds in an arbitrary braid system).

Let $\cal{F}^*$
be the union of systems $\cal{F}_n$. By construction, $\cal{F}^*$ is a complete
braid system. Moreover, $\cal{F}^*$ is the minimal braid system containing
$\cal{F}$.

Let $G=(\cal{A}, \{\phi,\e,\k,\sigma\})$ be a braided quantum group. According
to \eqref{21}--\eqref{22}, \eqref{235}--\eqref{236} and
Proposition~\ref{pro:28}
operators $\{\sigma,\tau\}$ form a braid system over the algebra $\cal{A}$.
The corresponding completion $\cal{F}=\{\sigma,\tau\}^*$ consists of maps
$\sigma_n\colon\cal{A}\otimes\cal{A}\rightarrow\cal{A}\otimes\cal{A}$ of the
form
\begin{equation}\label{A6}
\sigma_n=(\sigma\tau^{-1})^{n-1}\sigma=\sigma(\tau^{-1}\sigma)^{n-1}
\end{equation}
where $n\in\Bbb{Z}$.
\begin{pro}\label{pro:A2} The following identities hold
\begin{align}
(\phi\otimes\id)\sigma_{n+k}&=(\id\otimes\sigma_k)(\sigma_n\otimes\id)(
\id\otimes\phi)\label{A10}\\
\sigma_n(\id\otimes\k)&=(\k\otimes\id)\sigma_{-n}\label{A8}\\
\sigma_n(\k\otimes\id)&=(\id\otimes\k)\sigma_{-n}\label{A9}\\
(\id\otimes\phi)\sigma_{n+k}&=(\sigma_k\otimes\id)(\id\otimes\sigma_n)(
\phi\otimes\id).\label{A11}
\end{align}
\end{pro}
\begin{pf} Applying Proposition~\ref{pro:27} and \eqref{A6} we obtain
\begin{multline*}
\sigma_n(\id\otimes\k)=(\sigma\tau^{-1})^{n-1}\sigma(\id\otimes\k)=
(\k\otimes\id)(\tau\sigma^{-1})^{n-1}\tau\sigma^{-1}\tau\\
=(\k\otimes\id)(\sigma\tau^{-1})^{-n-1}\sigma=(\k\otimes\id)\sigma_{-n}.
\end{multline*}
Similarly,
$$\sigma_n(\k\otimes\id)=(\id\otimes\k)(\tau\sigma^{-1})^{n-1}\tau\sigma^{-1}
\tau=(\id\otimes\k)\sigma_{-n}.$$
Equalities \eqref{A10} and \eqref{A11} directly follow from
\eqref{220}--\eqref{223} and \eqref{229}--\eqref{230}. Indeed,
\begin{equation*}
\begin{split}
(\sigma_k\otimes\id)(\id\otimes\sigma_n)(\phi\otimes\id)&=
\bigl((\sigma\tau^{-1})^k\tau\otimes\id\bigr)\bigl(\id\otimes\sigma(\tau^{-1}
\sigma)^{n-1}\bigr)(\phi\otimes\id)\\
&=\bigl((\sigma\tau^{-1})^k\tau\otimes\id\bigr)(\id\otimes\sigma)
(\phi\otimes\id)(\tau^{-1}\sigma)^{n-1}\\
&=\bigl((\sigma\tau^{-1})^k\otimes\id\bigr)(\id\otimes\phi)\sigma
(\tau^{-1}\sigma)^{n-1}\\
&=(\id\otimes\phi)(\sigma\tau^{-1})^{n+k-1}\sigma=
(\id\otimes\phi)\sigma_{n+k}.
\end{split}
\end{equation*}
Similarly,
\begin{multline*}
(\id\otimes\sigma_n)(\sigma_k\otimes\id)(\id\otimes\phi)=
\bigl(\id\otimes(\sigma\tau^{-1})^n\tau\bigr)\bigl(\sigma(\tau^{-1}\sigma
)^{k-1}\otimes\id\bigr)(\id\otimes\phi)\\
=(\phi\otimes\id)(\sigma\tau^{-1})^{n+k-1}\sigma=(\phi\otimes
\id)\sigma_{n+k}.\qed
\end{multline*}
\renewcommand{\qed}{}
\end{pf}
As we shall now see, an arbitrary $\sigma_n\in\cal{F}$ is interpretable as
the flip-over operator corresponding to a modified braided
quantum group structure.

For each $n\in\Bbb{Z}$, let $m_n\colon\cal{A}\otimes\cal{A}\rightarrow
\cal{A}$ and $\k_n\colon\cal{A}\rightarrow\cal{A}$ be the maps given by
\begin{align}
m_n&=m\sigma_n^{-1}\sigma\label{A12}\\
\k_n&=(\e\otimes\k)\sigma_n^{-1}\sigma\phi=(\k\otimes\e)
\sigma_n^{-1}\sigma\phi \label{A13}
\end{align}
(the second equality in \eqref{A13} will be justified in the proof
of the proposition below). It is easy to see that each $m_n$, interpreted as
a product, determines a structure of an associative algebra on the space
$\cal{A}$. Indeed,
\begin{equation*}
\begin{split}
m_n(m_n\otimes\id)&=m\sigma_n^{-1}\sigma(m\sigma_n^{-1}\sigma\otimes\id)\\
&=m\sigma_n^{-1}(\id\otimes m)(\sigma\otimes\id)(\id\otimes\sigma)(
\sigma_n^{-1}\sigma\otimes\id)\\
&=m(m\otimes\id)(\id\otimes\sigma_n^{-1})(\sigma_n^{-1}\sigma\otimes\id)
(\id\otimes\sigma)(\sigma_n^{-1}\sigma\otimes\id)\\
&=m(m\otimes\id)(\id\otimes\sigma_n^{-1})(\sigma_n^{-1}\otimes\id)
(\id\otimes\sigma^{-1}_n)(\sigma\otimes\id)(\id\otimes\sigma)(\sigma\otimes
\id)\\
&=m(\id\otimes m)(\sigma_n^{-1}\otimes\id)(\id\otimes\sigma_n^{-1})
(\sigma_n^{-1}\otimes\id)(\id\otimes\sigma)(\sigma\otimes\id)
(\id\otimes\sigma)\\
&=m(\id\otimes m)(\sigma_n^{-1}\otimes\id)(\id\otimes\sigma_n^{-1}\sigma)
(\sigma\otimes\id)(\id\otimes\sigma_n^{-1}\sigma)\\
&=m\sigma_n^{-1}(m\otimes\id)(\id\otimes\sigma)(\sigma\otimes\id)
(\id\otimes\sigma_n^{-1}\sigma)\\
&=m\sigma_n^{-1}\sigma(\id\otimes m\sigma_n^{-1}\sigma)=m_n(\id\otimes m_n).
\end{split}
\end{equation*}

For each $n\in\Bbb{Z}$, let us denote by $\cal{A}_n$ the vector space
$\cal{A}$ endowed with the product $m_n$. Evidently, $1\in\cal{A}_n$ is
the unit in this algebra, too.
\begin{pro}\label{pro:A3} The pair $G_n=\bigl(\cal{A}_n,\{\phi,\e,
\k_n,\sigma_n\}\bigr)$ is a braided quantum group.
\end{pro}
\begin{pf} We have to check the last three axioms in
Definition~\ref{def:21}. The compatibility condition between
$\phi$ and $\sigma_n$ easily follows from \eqref{A10} and \eqref{A11}.
Further, a direct computation gives
\begin{equation*}
\begin{split}
\phi m_n&=(m\otimes m)(\id\otimes\sigma\otimes\id)(\phi\otimes
\phi)\sigma_n^{-1}\sigma=(m\otimes m)(\id\otimes\sigma\sigma_n^{-1}\sigma
\otimes\id)(\phi\otimes\phi)\\
&=(m\otimes m)(\id\otimes\sigma_{2-n}\otimes\id)(\phi\otimes\phi)\\
&=(m\sigma_n^{-1}\otimes m)(\id\otimes\phi\otimes\id)(\sigma_2\otimes\id)
(\id\otimes\phi)\\
&=(m\sigma_n^{-1}\sigma\otimes m)(\id\otimes\sigma\otimes\id)(\phi\otimes
\phi)\\
&=(m\sigma_n^{-1}\sigma\otimes m\sigma_n^{-1})(\id\otimes\phi\otimes\id)
(\id\otimes\sigma_{n+1})(\phi\otimes\id)\\
&=(m_n\otimes m_n)(\id\otimes\sigma_n\otimes\id)(\phi\otimes\phi).
\end{split}
\end{equation*}
Finally, we have to check that $k_n$ satisfies the antipode axiom. Let us
consider maps $k_n^\pm\colon\cal{A}\rightarrow\cal{A}$ given by
$$
k_n^-=(\k\otimes\e)\sigma_n^{-1}\sigma\phi\quad\qquad
k_n^+=(\e\otimes\k)\sigma_n^{-1}\sigma\phi.
$$
We have
\begin{equation*}
\begin{split}
m_n(k_n^-\otimes\id)\phi&=m\sigma_n^{-1}\sigma(\e\otimes\k\otimes\id)
(\tau\sigma_n^{-1}\sigma\otimes\id)(\phi\otimes\id)\phi\\
&=m(\e\otimes\k\otimes\id)(\id\otimes\sigma_{-n}^{-1}\sigma_{-1})
(\sigma_{1-n}\otimes\id)(\id\otimes\phi)\phi\\
&=m(\e\otimes\k\otimes\id)(\id\otimes\sigma_{-n}^{-1})(\phi\otimes
\id)\sigma_{-n}\phi\\
&=m(\e\otimes\k\otimes\id)(\tau\otimes\id)(\id\otimes\phi)\phi=
m(\k\otimes\id)\phi=1\e.
\end{split}
\end{equation*}
Similarly, it follows that
$m_n(\id\otimes\k_n^+)\phi=1\e. $
To complete the proof, let us observe that
\begin{multline*}
\k_n^+=(\e\otimes\k_n^+)\phi=m_n(m_n\otimes\id)(\k_n^-\otimes\id\otimes\k_n^+)
(\phi\otimes\id)\phi\\
=m_n(\id\otimes m_n)(\k_n^-\otimes\id\otimes\k_n^+)(\id\otimes\phi)\phi
=m_n(\k_n^-\otimes 1\e)\phi=k_n^-.
\end{multline*}
The map $k_n=k_n^\pm$ is bijective. Its inverse is given by
$$\k_n^{-1}\k=(\e\otimes\id)\sigma^{-1}\sigma_n\phi=(\id\otimes\e)\sigma^{-1}
\sigma_n\phi.\qed$$
\renewcommand{\qed}{}
\end{pf}

{}From the point of view of this analysis, the group $G_0$ is particularly
interesting. For example, left-covariant first-order differential structures
over $G$ (braided counterparts of structures considered in
\cite{WW}) are in a natural bijection with right $\cal{A}_0$-ideals $\cal{R}
\subseteq\ker(\e)$. Informally speaking, $G_0$ is interpretable
as a ``maximal braided simplification'' of $G$, with the same
coalgebra structure.
\section{A Special Case}
If $G$ is a standard quantum group then the counit map is geometrically
interpretable as a {\it classical point} of $G$. Let us analyze this question
in more details, in the braided context.
\begin{lem}
The following properties are equivalent
\begin{align}
\e m&=\e\otimes\e\label{M1}\\
(\e\otimes\id)\sigma&=\id\otimes\e\label{ML}\\
(\id\otimes\e)\sigma&=\e\otimes\id\label{MR}\\
\sigma&=\tau.\label{M3}
\end{align}
\end{lem}
\begin{pf}
Equality \eqref{M3} implies \eqref{M1},
according to \eqref{228}. If \eqref{M1} holds then \eqref{213}--\eqref{214}
imply \eqref{ML}--\eqref{MR}. Finally, if \eqref{ML} (or
\eqref{MR}) holds \eqref{216} implies that
two flip-over operators coincide.
\end{pf}

In other words, the mentioned conditions characterize the theory of
\cite{Maj}. The group $G_0$ introduced in the previous appendix is
of this kind (the braiding is given by $\tau$).

\end{document}